\documentclass[a4paper, 12pt]{article}
\usepackage{latexsym,amsmath,amssymb,bm, color}
\usepackage{amsmath,amssymb,color}

\newtheorem{thm}{Theorem}
\newtheorem{prop}[thm]{Proposition}
\newtheorem{lem}[thm]{Lemma}
\newtheorem{defi}[thm]{Definition}

\newcommand{\qed}{\hbox{\rule[-2pt]{3pt}{6pt}}}

\newcommand{\pf}{{\bf Proof. \ }}
\renewcommand{\qed}{\hfill $\Box$

\medskip}

\begin{document}
\title{\bf Matroids and Codes with the Rank Metric} 
\author{Keisuke Shiromoto\footnote{This work was supported by JSPS KAKENHI Grant Number 17K05348.}\\
Department of Mathematics and Engineering, \\
Kumamoto University, \\
2-39-1, Kurokami, Kumamoto 860-8555, Japan\\
           \tt{keisuke@kumamoto-u.ac.jp}
           } 
\maketitle
\begin{abstract}
We study the relationship between a $q$-analogue of matroids and linear codes with the rank metric in the vector space of matrices with entries in a finite field. We prove a Greene type identity for the rank generating function of these matroidal structures and the rank weight enumerator of these linear codes.
As an application, we give a combinatorial proof of a MacWilliams type identity for Delsarte rank-metric codes.\\



\end{abstract}

\section{Introduction}
Let $S$ be a finite set and let $\rho \: : \: 2^{S}\rightarrow {\mathbb Z}$ be a function.
An ordered pair $M=(S, \rho)$ is called a {\em matroid} if
\begin{enumerate}
\item[{\rm (1)}] If $X \subseteq S$, then $0 \leq \rho(X) \leq |X|$.
\item[{\rm (2)}] If $X,Y\subseteq S$ and $X \subseteq Y$, then $\rho(X) \leq \rho(Y)$.
\item[{\rm (3)}] If $X,Y\subseteq S$, then $\rho(X \cup Y)+\rho(X \cap Y) \leq \rho(X)+\rho(Y).
$
\end{enumerate}
Matroids are constructed from a lot of combinatorial or algebraic structures including graphs and matrices (cf. \cite{welsh76, oxley92}).
Let ${\mathbb F}_q$ be a finite filed of $q$ elements.
For an $[n,k]$ code $C$ over ${\mathbb F}_q$ and any subset $X \subseteq S:=\{1,2,\ldots,n\}$, define
$$
\rho(X):=\dim C\setminus (S-X) (=\dim C-\dim C(S-X)),
$$
where, for any  $Y \subseteq S$, $C\setminus Y$ denotes the punctured code by $Y$ and $$C(Y):=\{{\bm x} \in C \: : \: \mbox{\rm supp}({\bm x})\subseteq Y\}.$$
Then we find that $M_C:=(S,\rho)$ satisfies the above conditions and so $M_C$ is a matroid.

The (Whitney) {\em rank generating function} of a matroid $M=(S, \rho)$ is defined by
$$
R(M;x,y):=\sum_{X \subseteq S}x^{\rho(S)-\rho(X)}y^{|X|-\rho(X)}.
$$
In 1976, Greene (\cite{greene76}) proved one celebrated result, known as Greene's identity, on the relationship between the Hamming weight enumerator $W_C(x,y)$ of an $[n,k]$ code $C$ over ${\mathbb F}_q$ and the rank generating function of the corresponding matroid $M_C$ as follows:
$$
W_C(x,y)=y^{n-\dim C}(x-y)^{\dim C}R(M_C;\frac{qy}{x-y}, \frac{x-y}{y}).
$$
As an application of this result, he gave a elegant combinatorial proof of the MacWilliams identity for the Hamming weight enumerator (cf. \cite{welsh76}).
After this paper, a number of paper in the link on matroids and codes have been published, for instance, generalizations of the identity to higher weight enumerators (\cite{barg97, britz10, simonis94}), and to some enumerators for other classes of codes (\cite{BSW15, CGB13, simonis98, vertigan04}). 
In particular, Britz (\cite{britz10}) proved the equivalence between the set of higher weight enumerators of a linear code over a finite filed and the rank generating function associated to the code.
In \cite{BS16} and \cite{KMS17}, the authors studied the critical problem, known as a classical problem in matroid theory, from coding theoretical approach.

Delsarte introduced a rank-metric code as a set of matrices of given size over a finite field in \cite{delsarte78}.
He mainly studied the rank-metric code based on his research on association schemes.
In \cite{gabidulin85a, gabidulin85b}, Gabidulin gave a different definition of rank-metric codes such as a set of vectors in a vector space over an extension field.
The rank-weight of a codeword is defined by the rank of the associated matrix.
Ravagnani proved that Gabidulin codes can be regarded as a special case of Delsarte codes and he compared the duality theories of these codes in \cite{ravagnani2016}.
Some useful applications of rank-metric codes are widely known, for instance, to space-time codes (\cite{TaSeCa98}), to network coding (\cite{SiKsKo08}), and to cryptography (\cite{GaPaTr91}).

In addition, some of the results in classical coding theory are generalized to rank-metric codes. 
For instance, a Singleton bound on the minimum rank distance of the codes are proven in \cite{delsarte78}.
Some constructions of the codes which attain the bound, called {\em maximum rank distance} (MRD) {\em codes}, are known (\cite{CGLR18, delsarte78, gabidulin85a, gabidulin85b, KsGa05}).
The MacWilliams type identitis for rank weight enumerators of rank-metric codes are proposed as several expressions: an algebraic expression (\cite{delsarte78}), identitis for rank weight enumerators of Gabidulin rank-metric codes closed to the original form (\cite{gadouleau2007, gadouleau2008}), and moments of rank distributions of Delsarte rank-metric codes (\cite{ravagnani2016}).

The layout of this paper is as follows.
In Section 2, we firstly introduce some basic notion on Delsarte rank-metric codes and new matroidal structures related rank-metric codes.
And we give a Greene type identity for rank-metric codes.
Section 3 contains an application of the Greene type identity to the MacWilliams identity for Delsarte rank-metric codes.

Most of the basic terminology in coding theory and matroid theory used in this paper will be found in standard texts (cf. \cite{huffman2003, MS77, oxley92, welsh76}).

\section{$(q, r)$-Polymatroids and Rank Generating Functions}

We denote by $\mbox{\rm Mat}(n\times m, {\mathbb F}_q)$ the ${\mathbb F}_q$-vector space of $n \times m$ matrices with entries in ${\mathbb F}_q$.
A {\em Delsarte rank-metric code} ${\cal C}$ of size $n \times m$ over ${\mathbb F}_q$ is an ${\mathbb F}_q$-linear subspace of $\mbox{\rm Mat}(n\times m, {\mathbb F}_q)$.
Throughout this paper, set $E:={\mathbb F}_q^n$.
For convenience, we shall write for $D \leq E$ if $D$ is a subspace of $E$.
For any subspace $J \leq E$, define
\begin{eqnarray*}
{\cal C}(J) &:=&\{M \in {\cal C}\: |\: \mbox{\rm col}(M) \subseteq J\},\\
J^{\perp}&:=&\{{\bm y} \in {\mathbb F}_q^n \: | \: {\bm x} \cdot {\bm y}=0,\: {}^{\forall}\!{\bm x } \in J\},
\end{eqnarray*}
where $\mbox{\rm col}(M)$ denotes the column space of $M$ over ${\mathbb F}_q$.

\begin{lem}
\label{lem:subspaces}
${\cal C}(J)$ is an ${\mathbb F}_q$-linear subspace of $\mbox{\rm Mat}(n\times m, {\mathbb F}_q)$.
\end{lem}
\pf
For any $M, N \in {\cal C}(J)$, let ${\bm u}_1,\ldots, {\bm u}_m$ and ${\bm v}_1,\ldots, {\bm v}_m$ be the column vectors of $M$ and $N$, respectively.
For any $\alpha, \beta \in {\mathbb F}_q$, take a vector ${\bm x} \in \mbox{\rm col}(\alpha M+\beta N)$.
Then ${\bm x}$ can be written as
\begin{eqnarray*}
{\bm x}&=&a_1(\alpha {\bm u}_1+\beta {\bm v}_1)+\cdots+a_m(\alpha {\bm u}_m+\beta {\bm v}_m)\\
&=&\alpha(a_1 {\bm u}_1+\cdots+a_m {\bm u}_m)+\beta(a_1 {\bm v}_1+\cdots+a_m {\bm v}_m),
\end{eqnarray*}
for some $a_1,\ldots, a_m \in {\mathbb F}_q$.
Since $\alpha(a_1 {\bm u}_1+\cdots+a_m {\bm u}_m) \in J$ and $\beta(a_1 {\bm v}_1+\cdots+a_m {\bm v}_m) \in J$, we have that ${\bm x} \in J$.
\qed

We denote by $\Sigma(E)$ the set of all subspaces of $E$ and denote by $\Sigma(J)$ the set of all subspaces of any subspace $J \in \Sigma(E)$.

Now we shall introduce a $(q,r)$-polymatroid as a $q$-analogue of $k$-polymatroids (cf. \cite{Oxley-Whittle93}).
\begin{defi}
{\rm 
A $(q, r)$-{\em polymatroid} is an ordered pair $P=(E, \rho)$ consisting of a vector space $E:={\mathbb F}_q^n$ and a function $\rho \: : \: \Sigma(E)\to {\mathbb Z}^{+}\cup \{0\}$ having the following properties:
\begin{enumerate}
\item[{\rm (R1)}] If $A \leq E$, then $0 \leq \rho(A) \leq r \dim A$.
\item[{\rm (R2)}] If $A, B \leq E$ and $A \subseteq B$, then $\rho(A) \leq \rho(B)$.
\item[{\rm (R3)}] If $A, B \leq E$, then $\rho(A+B)+\rho(A \cap B) \leq \rho(A)+\rho(B).
$
\end{enumerate}
}
\end{defi}

We find easily that $\rho(\{{\bm 0}\})=0$ from (R1) and a $(q,1)$-polymatroid is a $q$-analogue of matroids.

\begin{prop}
\label{code-q-matroid}
Let ${\cal C}$ be a Delsarte rank-metric code in $\mbox{\rm Mat}(n\times m, {\mathbb F}_q)$, and let $\rho\: : \: \Sigma(E) \to {\mathbb Z}^{+}\cup \{0\}$ be the function defined by 
$$
\rho(J):=\dim_{{\mathbb F}_q} {\cal C}-\dim_{{\mathbb F}_q} {\cal C}(J^{\perp}),
$$
for any subspace $J$ of $E$.
Then $P_{\cal C}:=(E, \rho)$ is a $(q, m)$-polymatroid.
\end{prop}
\pf
Clearly the function $\rho$ satisfies (R1) and (R2).
Now we prove that (R3) holds.
Let $J_1$ and $J_2$ are subspaces of $E$.
From Lemma \ref{lem:subspaces}, we have that
$$
\dim {\cal C}(J_1^{\perp})+\dim {\cal C}(J_2^{\perp})=\dim {\cal C}(J_1^{\perp})\cap {\cal C}(J_2^{\perp})+\dim({\cal C}(J_1^{\perp})+{\cal C}(J_2^{\perp})).
$$
Moreover, it follows that
\begin{eqnarray*}
&& {\cal C}(J_1^{\perp})\cap {\cal C}(J_2^{\perp})={\cal C}(J_1^{\perp} \cap J_2^{\perp})={\cal C}((J_1+J_2)^{\perp}),\\
&& {\cal C}(J_1^{\perp})+ {\cal C}(J_2^{\perp}) \subseteq {\cal C}(J_1^{\perp}+J_2^{\perp})={\cal C}((J_1 \cap J_2)^{\perp}),
\end{eqnarray*}
because of $\mbox{\rm col}(M_1+M_2) \subseteq \mbox{\rm col}(M_1)+\mbox{\rm col}(M_2)$ for any $M_1, M_2 \in \mbox{\rm Mat}(n\times m, {\mathbb F}_q)$.
Therefore we have that
\begin{eqnarray*}
\rho(J_1)+\rho(J_2)&=&(\dim_{{\mathbb F}_q} {\cal C}-\dim_{{\mathbb F}_q} {\cal C}(J_1^{\perp}))+(\dim_{{\mathbb F}_q} {\cal C}-\dim_{{\mathbb F}_q} {\cal C}(J_2^{\perp}))\\
&\geq&(\dim_{{\mathbb F}_q} {\cal C}-\dim_{{\mathbb F}_q} {\cal C}((J_1+J_2)^{\perp}))
+(\dim_{{\mathbb F}_q} {\cal C}-\dim_{{\mathbb F}_q} {\cal C}((J_1\cap J_2)^{\perp}))\\
&=&\rho(J_1+J_2)+\rho(J_1\cap J_2).
\end{eqnarray*}
Thus $P_{{\cal C}}=(E, \rho)$ is a $(q, m)$-polymatroid.
\qed

The {\em rank generating function} of  a $(q, r)$-polymatroid $P=(E, \rho)$ is defined by
$$
R_{P}(X_1, X_2, X_3, X_4):=\sum_{D \in \Sigma(E)}f_{P}^D(X_1, X_2)g^{\dim D}(X_3, X_4),
$$
where
$$
f_{P}^J(X, Y):=X^{\rho(E)-\rho(J)}Y^{r\dim J-\rho(J)},
$$
for any subspace $J \in \Sigma(E)$, and
$$
g^{l}(X, Y):=\prod_{i=0}^{l-1}(X-q^iY),
$$
for any nonnegative integer $l \in {\mathbb Z}^{+}\cup\{0\}$.

The following lemma is essential.
\begin{lem}
\label{lem:rank}
Let $P=(E,\rho)$ be a $(q, r)$-polymatroid.
If $A$ and $B$ are subspaces of $E$ such as $A \subseteq B$, then
$$
0 \leq \rho(A)-\rho(B)-r(\dim A-\dim B).
$$
\end{lem}
\pf
Let $D$ be a subspace of $E$ generated by $B-A$.
From (R1), we have that
$$
\rho(D) \leq r\dim D=r(\dim B-\dim A).
$$
From (R3), we also have that
$$
\rho(A)+\rho(D)\geq \rho(A+D)+\rho(A\cap D)\geq\rho(B).
$$
By combining the above two inequalities, the lemma follows.
\qed

\begin{prop}
\label{dual}
For any $(q, r)$-polymatroid $P=(E, \rho)$ and any subspace $J$ of $E$, define
$$
\rho^*(J):=\rho(J^{\perp})+r \dim J-\rho(E).
$$
Then $P^*=(E, \rho^*)$ is a $(q, r)$-polymatroid.
\end{prop}
\pf
$J^{\perp} \subseteq E$ implies that $\rho(J^{\perp}) \leq \rho(E)$.
So we have that $(0\leq)\rho^*(J)\leq r \dim J$.
Thus (R1) holds.

Let $J_1$ and $J_2$ are subspaces of $E$.
Assume that $J_1 \subseteq J_2$. 
Then it follows from Lemma \ref{lem:rank} that
\begin{eqnarray*}
\rho^*(J_2)-\rho^*(J_1)&=&(\rho(J_2^{\perp})+r \dim J_2-\rho(E))-(\rho(J_1^{\perp})+r \dim J_1-\rho(E))\\
&=&\rho(J_2^{\perp})-\rho(J_1^{\perp})+r(\dim J_2-\dim J_1)\\
&=&\rho(J_2^{\perp})-\rho(J_1^{\perp})+r(\dim J_1^{\perp}-\dim J_2^{\perp})\\
&\geq&0.
\end{eqnarray*}
Hence (R2) holds.

From the definition of $\rho^*$, it follows that
\begin{eqnarray*}
\rho^*(J_1\cap J_2)&=&\rho((J_1\cap J_2)^{\perp})+r\dim (J_1\cap J_2)-\rho(E),\\
\rho^*(J_1+ J_2)&=&\rho((J_1+ J_2)^{\perp})+r\dim (J_1+ J_2)-\rho(E),\\
\rho^*(J_1)&=&\rho(J_1^{\perp})+r\dim J_1-\rho(E),\\
\rho^*(J_2)&=&\rho(J_2^{\perp})+r\dim J_2-\rho(E).
\end{eqnarray*}
From the inequality in (R3) for $\rho$, we have that
\begin{eqnarray*}
&&\rho^*(J_1\cap J_2)+\rho^*(J_1+ J_2)-\rho^*(J_1)-\rho^*(J_2)\\
&=&\{\rho(J_1^{\perp}+ J_2^{\perp})+\rho(J_1^{\perp}\cap J_2^{\perp})
 -(\rho(J_1^{\perp})+\rho(J_2^{\perp}))\}\\
 &&\hfill+r\{\dim (J_1\cap J_2)+\dim (J_1+ J_2)-(\dim J_1+\dim J_2)\}\\
&\leq& 0.
\end{eqnarray*}
Then (R3) holds.
\qed

For a rank generating function $R_{P}(X_1, X_2, X_3, X_4)$ of  a $(q, r)$-polymatroid $P=(E, \rho)$, define
$$
\widehat{R}_{P}(X_1, X_2, X_3, X_4):=\sum_{D \in \Sigma(E)}f_{P}^D(X_1, X_2)g^{\dim D^{\perp}}(X_3, X_4).
$$
The following result is a duality on rank generating functions of $(q, r)$-polymatroids.

\begin{thm}
\label{thm:duality}
Let $P=(E, \rho)$ be a $(q, r)$-polymatroid.
Then
$$
R_{P^*}(X_1, X_2, X_3, X_4)=\widehat{R}_{P}(X_2, X_1, X_3, X_4).
$$
\end{thm}
\pf
It follows from the definition of rank generating functions that
\begin{eqnarray*}
R_{P^*}(X_1, X_2, X_3, X_4)&=&\sum_{D \in \Sigma(E)}f_{P^*}^D(X_1, X_2)g^{\dim D}(X_3, X_4)\\
&=&\sum_{D \in \Sigma(E)}X_1^{\rho^*(E)-\rho^*(D)}X_2^{r\dim D-\rho^*(D)}
\prod_{i=0}^{\dim D-1}(X_3-q^iX_4)\\
&=&\sum_{D \in \Sigma(E)}X_1^{(rn-\rho(E))-(\rho(D^{\perp})+r\dim D-\rho(E))}X_2^{r\dim D-(\rho(D^{\perp})+r\dim D-\rho(E))}\\
&& \hspace{8cm} \times \prod_{i=0}^{\dim D-1}(X_3-q^iX_4)\\
&=&\sum_{D \in \Sigma(E)}X_1^{r(n-\dim D)-\rho(D^{\perp})}X_2^{\rho(E)-\rho(D^{\perp})}
\prod_{i=0}^{\dim D-1}(X_3-q^iX_4)\\
&=&\sum_{D' \in \Sigma(E)}X_2^{\rho(E)-\rho(D')}X_1^{r\dim D'-\rho(D')}
\prod_{i=0}^{n-\dim D'-1}(X_3-q^iX_4)\\
&=&\widehat{R}_{P}(X_2, X_1, X_3, X_4).
\end{eqnarray*}
\qed

The {\em trace product} of matrices $M, N \in \mbox{\rm Mat}(n\times m, {\mathbb F}_q)$
is defined by
$$
\langle M, N\rangle :=\mbox{\rm Tr}(MN^{{\mbox{\rm T}}}),
$$
where $\mbox{\rm Tr}$ denotes the trace of a matrix.
Then it is easy to find that the trace product can be calculate from the standard inner products of each row vectors as follows:

\begin{lem}
\label{lem:trace-products}
Let $M_i$ and $N_i$ be the $i$-th column vectors of $M$ and $N$, respectively.
Then we have that
$$
\langle M, N\rangle=\sum_{i=1}^m M_i\cdot N_i,
$$
where $M_i\cdot N_i$ denotes the standard inner product of vectors $M_i$ and $N_i$.
\end{lem}

Let ${\cal C}$ be a Delsarte rank-metric code in $\mbox{\rm Mat}(n\times m, {\mathbb F}_q)$.
The {\em dual code} of ${\cal C}$ is defined by
$$
{\cal C}^{\perp}:=\{N \in \mbox{\rm Mat}(n\times m, {\mathbb F}_q)\: : \: \mbox{\rm $\langle M, N\rangle=0$ for all $M \in {\cal C}$}\}.
$$
Then the following result is well-known (see, Lemma 5 in \cite{ravagnani2016}).

\begin{lem} 
\label{lem:dual-codes}
{\rm (\cite{ravagnani2016})}
For any Delsarte rank-metric codes ${\cal C}, {\cal D} \in \mbox{\rm Mat}(n\times m, {\mathbb F}_q)$,
it follows that
\begin{enumerate}
\item[{\rm (1)}] $({\cal C}^{\perp})^{\perp}={\cal C}$;
\item[{\rm (2)}] $\dim _{{\mathbb F}_q} ({\cal C}^{\perp})=nm-\dim _{{\mathbb F}_q} ({\cal C})$;
\item[{\rm (3)}] $({\cal C} \cap {\cal D})^{\perp}={\cal C}^{\perp}+{\cal D}^{\perp}$, and $({\cal C} + {\cal D})^{\perp}={\cal C}^{\perp}\cap{\cal D}^{\perp}$.
\end{enumerate}
\end{lem}

Set ${\cal V}:=\mbox{\rm Mat}(n\times m, {\mathbb F}_q)$.
Then it turns out easily that, for any subspace $R$ of ${\mathbb F}_q^n$,
\begin{enumerate}
\item[(1)] ${\cal C}(R)={\cal V}(R)\cap {\cal C}$, for any Delsarte rank-metric code ${\cal C} \subseteq {\cal V}$;
\item[(2)] $|{\cal V}(R)|=|R|^m=q^{m\dim R}$.
\end{enumerate}

\begin{lem}
\label{lem: duality}
For any subspace $R$ of ${\mathbb F}_q^n$, 
$$
({\cal V}(R))^{\perp}={\cal V}(R^{\perp}).
$$
\end{lem}
\pf
For any matrix $M\in {\cal V}(R^{\perp})$, Lemma \ref{lem:trace-products} implies that
\begin{eqnarray*}
\langle M, N\rangle &=&\mbox{\rm Tr}(MN^{{\mbox{\rm T}}})\\
&=&\sum_{i=1}^m M_i\cdot N_i=0,
\end{eqnarray*}
for all $N \in {\cal V}(R)$.
Thus we have that $M \in ({\cal V}(R))^{\perp}$ and so 
${\cal V}(R^{\perp}) \subseteq ({\cal V}(R))^{\perp}$.

Conversely, take a matrix $M \in  ({\cal V}(R))^{\perp}$.
Set
$$
R=\{{\bm x}_1,\ldots,{\bm x}_{|R|}\}.
$$
Let $M({\bm x}_i, j)$ be the $n \times m$ matrix over ${\mathbb F}_q$ with ${\bm x}_i$ as the $j$-th column vector and ${\bm 0}$ elsewhere.
Then it turns out that $M({\bm x}_i, j) \in {\cal V}(R)$.
Thus we have that
\begin{eqnarray*}
0&=& \langle M, M({\bm x}_i, j)\rangle=\sum_{l=1}^m M_l\cdot (M({\bm x}_i, j))_l\\
&=&M_j \cdot {\bm x}_i.
\end{eqnarray*}
This implies that $M_j \in R^{\perp}$, for all $j=1,2,\ldots, m$.
Therefore it follows that $M \in {\cal V}(R^{\perp})$ and so ${\cal V}(R^{\perp}) \supseteq ({\cal V}(R))^{\perp}$.
\qed

Let ${\cal C}$ be a Delsarte rank-metric code in $\mbox{\rm Mat}(n\times m, {\mathbb F}_q)$.
The {\em dual space} ${\cal C}^{\star}$ of ${\cal C}$ is defined by
$$
{\cal C}^{\star}:=\mbox{\rm Hom}_{{\mathbb F}_q}({\cal C}, {\mathbb F}_q).
$$
Since ${\cal C}$ is a finite free ${\mathbb F}_q$-module, there is a (non-natural) isomorphism ${\cal C} \cong {\cal C}^{\star}$ (see, for instance, \cite{lam99}).

The following exact sequence is an analogue of the basic exact sequence on Proposition 2 in \cite{shiromoto99}.
 
\begin{prop}
\label{thm:basic-exact-sequence}
Let ${\cal C}$ be a Delsarte rank-metric code in $\mbox{\rm Mat}(n\times m, {\mathbb F}_q)$.
For any subspace $R$ of ${\mathbb F}_q^n$, there is an exact sequence as ${\mathbb F}_q$-modules,
$$
0 \longrightarrow {\cal C}^{\perp}(R) \xrightarrow{\scriptsize {\mbox{\rm inc}}} {\cal V}(R) \xrightarrow{ f }
{\cal C}^{\star} \xrightarrow{\scriptsize {\mbox{\rm res}}} {\cal C}(R^{\perp})^{\star} \longrightarrow 0,
$$
where $f$ is an ${\mathbb F}_q$-homomorphism defined by
\begin{eqnarray*}
f\: : \: {\cal V} &\longrightarrow& {\cal C}^{\star}\\
M &\longmapsto& (\hat{M}\: : \: N \longmapsto \langle M, N\rangle),
\end{eqnarray*}
and the maps {\rm inc} and {\rm res} denote the inclusion map and the restriction map, respectively.
\end{prop}
\pf
The inclusion map {\rm inc} is a natural injection, and the restriction map {\rm res} is surjective, since 
${\mathbb F}_q$ is an injective module over itself.

Take $M \in {\cal C}^{\perp}(R)$.
Then it follows that
$$
\hat{M}\: : \: N \longmapsto \langle M, N\rangle=0,
$$
for all $N \in {\cal C}$.
This implies that $M \in \ker f$ and so $\mbox{\rm Im}(\mbox{\rm inc})\subseteq \ker f$.
Conversely, take $M \in \ker f$.
Then we have that $\hat{M}(N)=\langle M, N\rangle=0$ for all $N \in {\cal C}$ and so $M \in {\cal C}^{\perp} \cap {\cal V}(R)={\cal C}^{\perp}(R)$.
Therefore we have that $\mbox{\rm Im}(\mbox{\rm inc})\supseteq \ker f$.
Thus it follows that $\mbox{\rm Im}(\mbox{\rm inc})= \ker f$.

Next we prove that $\ker (\mbox{\rm res})= \mbox{\rm Im} f $.
Take $M \in {\cal V}(R)$.
Then it also follows from Lemma \ref{lem:trace-products} that $\hat{M}(N)=\langle M, N\rangle=0$ for all $N \in {\cal C}(R^{\perp})$.
This implies that $\hat{M} \in \ker (\mbox{\rm res})$ and so $\mbox{\rm Im} f \subseteq \ker (\mbox{\rm res})$.
Conversely, take $\lambda \in \ker (\mbox{\rm res})$.
The map $f$ is surjective.
Thus there exists $M \in {\cal V}$ such that $\lambda=\hat{M}$.
Therefore we have that $\hat{M}(N)=\langle M, N\rangle=0$ for all $N \in {\cal C}(R^{\perp})$, and so $M \in ({\cal C}(R^{\perp}))^{\perp}$.
From Lemma \ref{lem: duality}, we have that
\begin{eqnarray*}
({\cal C}(R^{\perp}))^{\perp}&=&({\cal C}\cap {\cal V}(R^{\perp}))^{\perp}={\cal C}^{\perp}+{\cal V}(R^{\perp})^{\perp}\\
&=&{\cal C}^{\perp}+{\cal V}(R),
\end{eqnarray*}
and so $M \in {\cal C}^{\perp}+{\cal V}(R)$.
Then $M=M_1+M_2$ for some $M_1 \in {\cal C}^{\perp}$ and $M_2 \in {\cal V}(R)$.
Since $\hat{M_1}(N)=0$ for all $N \in {\cal C}$, it turns out that 
$\lambda=\hat{M}=\hat{M_2} \in \mbox{\rm Im} f$.
Thus it follows that $\mbox{\rm Im} f= \ker (\mbox{\rm res})$.

Therefore the proposition follows.
\qed

\begin{prop}
\label{prop:duality}
For any Delsarte rank-metric code ${\cal C}$ in $\mbox{\rm Mat}(n\times m, {\mathbb F}_q)$, let $P_{\cal C}=(E, \rho)$ be the $(q, m)$-polymatroid obtained from ${\cal C}$ such as discussed in Proposition {\rm \ref{code-q-matroid}}.
Then $P_{\cal C}^{*}=P_{{\cal C}^{\perp}}$.
\end{prop}
\pf
Take any subspace $J \subseteq E$.
Set
$$
\tau(J):=\dim_{{\mathbb F}_q} {\cal C}^{\perp} -\dim_{{\mathbb F}_q} {\cal C}^{\perp}(J^{\perp}).
$$
From the definition of dual $(q, m)$-polymatroids, we have that
\begin{eqnarray*}
\rho^{*}(J)&=&\rho(J^{\perp})+m \dim J-\rho(E)\\
&=&(\dim_{{\mathbb F}_q} {\cal C}-\dim_{{\mathbb F}_q} {\cal C}(J))+m \dim J-\dim_{{\mathbb F}_q} {\cal C}\\
&=& m \dim J -\dim_{{\mathbb F}_q} {\cal C}(J).
\end{eqnarray*}
By applying $J^{\perp}$ for $R$ in Proposition \ref{thm:basic-exact-sequence}, we have that
\begin{eqnarray*}
&&|{\cal C}^{\perp}(J^{\perp})|\cdot |{\cal C}^{\star}|=|{\cal V}(J^{\perp})|\cdot |{\cal C}(J)^{\star}|\\
&\Longleftrightarrow& |{\cal C}^{\perp}(J^{\perp})|\cdot |{\cal C}|=|{\cal V}(J^{\perp})|\cdot |{\cal C}(J)|\\
&\Longleftrightarrow& q^{\dim_{{\mathbb F}_q} {\cal C}^{\perp}(J^{\perp})} \times q^{\dim_{{\mathbb F}_q} {\cal C}}=q^{m \dim J^{\perp}} \times q^{\dim_{{\mathbb F}_q} {\cal C}(J)}\\
&\Longleftrightarrow& \dim_{{\mathbb F}_q} {\cal C}^{\perp}(J^{\perp})+\dim_{{\mathbb F}_q} {\cal C}=m \dim J^{\perp}+\dim_{{\mathbb F}_q} {\cal C}(J)
\end{eqnarray*}
From the above equation and Lemma \ref{lem:dual-codes}(2), it follows that
\begin{eqnarray*}
\tau(J)&=& (nm-\dim_{{\mathbb F}_q} {\cal C})-(m \dim J^{\perp}+\dim_{{\mathbb F}_q} {\cal C}(J)-\dim_{{\mathbb F}_q} {\cal C})\\
&=& m \dim J-\dim_{{\mathbb F}_q} {\cal C}(J)\\
&=& \rho^{*}(J).
\end{eqnarray*}
\qed

\begin{defi}
{\rm (cf. \cite{gadouleau2008, ravagnani2016})
Let ${\cal C}$ be a Delsarte rank-metric code in $\mbox{\rm Mat}(n\times m, {\mathbb F}_q)$.
The {\em rank distribution} of ${\cal C}$ is the collection $\{A_i({\cal C})\}_{i\in {\mathbb Z}_{\geq 0}}$, where
$$
A_i({\cal C}):=|\{M \in {\cal C} \: : \: \mbox{\rm rank}(M)=i\}|.
$$
The {\em rank weight enumerator} of ${\cal C}$ is defined by
$$
W_{{\cal C}}^{{\rm R}} (x, y):=\sum_{i=0}^n A_i({\cal C}) x^{n-i}y^i.
$$
}
\end{defi}

For any Delsarte rank-metric code ${\cal C}$ in $\mbox{\rm Mat}(n\times m, {\mathbb F}_q)$ and any subspace $R \subseteq  {\mathbb F}_q^n$, set
\begin{eqnarray*}
A_{{\cal C}}(R)&:=&|\{M \in {\cal C} \: : \: \mbox{\rm col}(M)=R\}|,\\
B_{{\cal C}}(R)&:=&|\{M \in {\cal C} \: : \: \mbox{\rm col}(M)\subseteq R\}|.
\end{eqnarray*}
Then it turns out easily that $B_{{\cal C}}(R)=\sum_{T \in \Sigma(R)}A_{{\cal C}}(T)$ and $B_{{\cal C}}(R)=|{\cal C}(R)|$.
From the M\"obius inversion formula (see, for instance, \cite{aigner79}), we have the following equation.

\begin{lem}
\label{lem:Mobius-inversion}
$$
A_{{\cal C}}(R)=\sum_{T \in \Sigma(R)}(-1)^{\dim R-\dim T} q^{{\dim R - \dim T \choose 2}}B_{{\cal C}}(T),
$$
where ${a \choose b}$ denotes the binomial coefficient of integers $a \geq b \geq 0$.
\end{lem}

Now we prove our main result in this paper.
The following equation is a kind of Greene type identities for $(q, r)$-polymatroids and Delsarte rank-metric codes.

\begin{thm}
\label{thm:greene-type-id}
Let ${\cal C}$ be a Delsarte rank-metric code in $\mbox{\rm Mat}(n\times m, {\mathbb F}_q)$, and let $P_{\cal C}=(E, \rho)$ be the $(q, m)$-polymatroid obtained from ${\cal C}$ such as discussed in Proposition {\rm \ref{code-q-matroid}}.
Then 
$$
W_{{\cal C}}^{{\rm R}} (x, y)=y^{n-\dim {\cal C} / m}
R_{P_{\cal C}}(qy^{1/m}, \frac{1}{y^{1/m}}, x, y).
$$
\end{thm}
\pf
By Lemma \ref{lem:Mobius-inversion} and the $q$-binomial Theorem, we obtain that
\begin{eqnarray*}
W_{{\cal C}}^{{\rm R}} (x, y)&=& \sum_{i=0}^n A_i({\cal C}) x^{n-i}y^i\\
&=&\sum_{R \in \Sigma(E)}A_{{\cal C}}(R)x^{n-\dim R}y^{\dim R}\\
&=&\sum_{R \in \Sigma(E)}\left( \sum_{T \in \Sigma(R)}(-1)^{\dim R-\dim T} q^{{\dim R - \dim T \choose 2}}B_{{\cal C}}(T)\right)x^{n-\dim R}y^{\dim R}\\
&=&\sum_{T \in \Sigma(E)}B_{{\cal C}}(T)\left(\sum_{\substack{T \leq R \\ R \in \Sigma(E)}}(-1)^{\dim R-\dim T} q^{{\dim R - \dim T \choose 2}}\right)x^{n-\dim R}y^{\dim R}\\
&=&\sum_{T \in \Sigma(E)}B_{{\cal C}}(T)\sum_{J \leq T^{\perp}}(-1)^{\dim J} q^{{\dim J \choose 2}}x^{n-\dim J-\dim T}y^{\dim J+\dim T}\\
&=&\sum_{T \in \Sigma(E)}B_{{\cal C}}(T)\sum_{J \in \Sigma(T^{\perp})}(-1)^{\dim J} q^{{\dim J \choose 2}}x^{\dim T^{\perp}-\dim J}y^{\dim J+\dim T}\\
&=&\sum_{T \in \Sigma(E)}B_{{\cal C}}(T)y^{\dim T}\prod_{j=0}^{\dim T^{\perp}-1}(x-q^jy)\\
&=&\sum_{T' \in \Sigma(E)}B_{{\cal C}}((T')^{\perp})y^{n-\dim T'}\prod_{j=0}^{\dim T'-1}(x-q^jy)\\
&=&\sum_{T' \in \Sigma(E)}|{\cal C}((T')^{\perp})|y^{n-\dim T'}\prod_{j=0}^{\dim T'-1}(x-q^jy)\\
&=&\sum_{T \in \Sigma(E)}q^{\dim_{{\mathbb F}_q} {\cal C}(T^{\perp})}y^{n-\dim T}\prod_{j=0}^{\dim T-1}(x-q^jy)\\
&=&\sum_{T \in \Sigma(E)}q^{\dim_{{\mathbb F}_q} {\cal C}-\rho (T)}y^{n-\dim T}\prod_{j=0}^{\dim T-1}(x-q^jy)\\
&=&y^{n-\rho (E)/m}\sum_{T \in \Sigma(E)}\left(\frac{1}{y^{1/m}}\right)^{m \dim T-\rho (T)} \left( qy^{1/m}\right)^{\rho (E)-\rho(T)}\prod_{j=0}^{\dim T-1}(x-q^jy)\\
&=&y^{n-\dim_{{\mathbb F}_q} {\cal C}/m}\sum_{T \in \Sigma(E)}f_{P_{{\cal C}}}^T\left(qy^{1/m}, \frac{1}{y^{1/m}}\right)g^{\dim T}(x, y).
\end{eqnarray*}
\qed

\section{Application to the MacWilliams Identity}

In this section, we shall present a MacWilliams type identity for Delsarte rank-metric codes as an application of Theorem \ref{thm:greene-type-id}.

\begin{prop}
\label{prop:dual-enumerator}
Let ${\cal C}$ be a Delsarte rank-metric code in $\mbox{\rm Mat}(n\times m, {\mathbb F}_q)$.
Then 
$$
W_{{\cal C}^{\perp}}^{{\rm R}} (x, y)=\frac{1}{|{\cal C}|}
\sum_{S \in \Sigma(E)}A_{{\cal C}}(S)\sum_{j=0}^n \sum_{l=0}^{j} {n-\dim S \brack j-l}_q {n-j+l \brack l}_q  (-1)^l q^{{l \choose 2}} q^{m(j-l)}  y^j x^{n-j},
$$
where ${a \brack b}_q$ denote the Gaussian binomial coefficient (or $q$-binomial coefficient) of non-nagative integers $a$ and $b$.
\end{prop}
\pf
By combining Theorem \ref{thm:greene-type-id}, Proposition \ref{prop:duality} and Theorem \ref{thm:duality}, we have that
\begin{eqnarray*}
W_{{\cal C}^{\perp}}^{{\rm R}} (x, y)&=&y^{n-\dim_{{\mathbb F}_q} {\cal C}^{\perp} / m}
\sum_{D \in \Sigma(E)}f_{P_{{\cal C}^{\perp}}}^D\left(qy^{1/m}, \frac{1}{y^{1/m}}\right)g^{\dim D}(x, y)\\
&=&y^{\dim_{{\mathbb F}_q} {\cal C} / m}
\sum_{D \in \Sigma(E)}f_{P_{{\cal C}}}^D\left(\frac{1}{y^{1/m}}, qy^{1/m}\right)g^{n-\dim D}(x, y)\\
&=&y^{\dim_{{\mathbb F}_q} {\cal C} / m}
\sum_{D \in \Sigma(E)} \left( qy^{1/m}\right)^{m \dim D-\rho (D)} \left(\frac{1}{y^{1/m}}\right)^{\rho (E)-\rho(D)}\prod_{j=0}^{n- \dim D-1}(x-q^jy)\\
&=&q^{-\dim_{{\mathbb F}_q} {\cal C}}
\sum_{D \in \Sigma(E)} q^{\dim_{{\mathbb F}_q} {\cal C}-\rho(D)} \left( q^my\right)^{\dim D} \prod_{j=0}^{n- \dim D-1}(x-q^jy)\\
&=& \frac{1}{|{\cal C}|}
\sum_{D' \in \Sigma(E)} q^{\dim {\cal C}(D')} \left( q^my\right)^{n-\dim D'} \prod_{j=0}^{\dim D'-1}(x-q^jy)\\
&=&\frac{1}{|{\cal C}|}
\sum_{D \in \Sigma(E)}B_{{\cal C}}(D) (q^my)^{n-\dim D} \prod_{j=0}^{\dim D-1}(x-q^jy)
\\
&=&\frac{1}{|{\cal C}|}
\sum_{D \in \Sigma(E)}\left(\sum_{S \in \Sigma(D)} A_{{\cal C}}(S)\right) (q^my)^{n-\dim D} \sum_{u=0}^{\dim D} {\dim D \brack u}_q (-1)^u q^{{u \choose 2}}y^u x^{\dim D-u}\\
&=& \frac{1}{|{\cal C}|}
\sum_{S \in \Sigma(E)}A_{{\cal C}}(S)\sum_{j=0}^n {n-\dim S \brack n-j}_q q^{m(n-j)} \sum_{l=0}^{j} {j \brack l}_q (-1)^l q^{{l \choose 2}}y^{n-j+l} x^{j-l}\\
&=& \frac{1}{|{\cal C}|}
\sum_{S \in \Sigma(E)}A_{{\cal C}}(S)\sum_{j=0}^n \sum_{l=0}^{j} {n-\dim S \brack j-l}_q {n-j+l \brack l}_q  (-1)^l q^{{l \choose 2}} q^{m(j-l)}  y^j x^{n-j}.
\end{eqnarray*}
\qed

The equation in Proposition \ref{prop:dual-enumerator} can be viewed as a kind of MacWilliams type identities for Delsarte rank-matric codes.
Now we prove the equivalence between the above equation and the MacWilliams type identity for Gabidulin rank-metric codes proposed in \cite{gadouleau2007, gadouleau2008} by using the concepts of $q$-products and $q$-derivative for homogeneous polynomials.
 
\begin{defi}
\label{defi:q-product}
{\rm
(\cite{gadouleau2007, gadouleau2008})
Let $a(x,y;m)=\sum_{i=0}^r a_i(m)x^{r-i}y^i$ and $b(x,y;m)=\sum_{j=0}^s b_j(m)x^{s-j}y^j$ be two homogeneous polynomials in $x$ and $y$ of degree $r$ and $s$ respectively with coefficients $a_i(m)$ and $b_j(m)$ for $i, j \geq 0$ in turn are real functions of $m$, and are assumed to be zero unless otherwise specified.
The $q$-{\em product} $c(x,y;m)$ of $a(x,y;m)$ and $b(x,y;m)$ is defined to be the homogeneous polynomial of degree $(r+s)$ as follows:
$$
c(x,y;m):=a(x,y;m)*b(x,y;m)=\sum_{u=0}^{r+s}c_u(m)x^{r+s-u}y^u,
$$
where
$$
c_u(m):=\sum_{i=0}^uq^{is}a_i(m)b_{u-i}(m-i).
$$
For $n\geq 0$, the {\em $n$-th $q$-power} of $a(x,y;m)$ is defined recursively:
\begin{eqnarray*}
a(x,y;m)^{[0]}&:=&1,\: \mbox{\rm and}\\
a(x,y;m)^{[n]}&:=&a(x,y;m)^{[n-1]}*a(x,y;m),\: n=1,2,\ldots,
\end{eqnarray*}
}
\end{defi}

\begin{lem}
\label{lem:q-products}
For any $n \geq 0$ and $(0 \leq) l \leq n$,
\begin{enumerate}
\item[{\rm (1)}] $(q^m y)^{[n]}=(q^my)^n$,
\item[{\rm (2)}] $(x-y)^{[l]}*(q^my)^{[n-l]}=(x-y)^{[l]}\times (q^my)^{n-l}$.
\end{enumerate}
\end{lem}
\pf
(1) We shall prove this equation by induction on $n$.
When $n=0$, it follows that $(q^m y)^{[0]}=1=(q^my)^0$.
Suppose that the equation is true for $n-1$, that is,
$(q^m y)^{[n-1]}=(q^my)^{n-1}$.
Then, from the definition of $q$-products, we have that
\begin{eqnarray*}
(q^m y)^{[n]}&=&(q^m y)^{[n-1]}*(q^m y)=(q^my)^{n-1}*q^m y\\
&=&(q^{n-1}\times q^{m(n-1)}\times q^{m-n+1})y^n=(q^my)^n.
\end{eqnarray*}
(2) From Lemma 1 in \cite{gadouleau2007}, it follows that
$$
(x-y)^{[l]}=\sum_{u=0}^l {l \brack u}_q (-1)^uq^{{u \choose 2}}x^{l-u}y^u,
$$
where
$$
{l \brack u}_q:=\prod_{i=0}^{u-1}\frac{1-q^{l-i}}{1-q^{i+1}}.
$$
Then we have that
\begin{eqnarray*}
(x-y)^{[l]}*(q^my)^{[n-l]}&=&\left(\sum_{u=0}^l {l \brack u}_q (-1)^uq^{{u \choose 2}}x^{l-u}y^u\right)*\left(q^{m(n-l)}y^{n-l}\right)\\
&=&\sum_{u=0}^n c_u(m) x^{n-u}y^u,
\end{eqnarray*}
where
\begin{eqnarray*}
c_u(m)&=&0,\: u=0,1,\ldots,n-l-1,\\
c_{n-l+j}(m)&=&q^{m(n-l)}{l \brack j}_q (-1)^jq^{{j \choose 2}},\: j=0,1,\ldots,l.
\end{eqnarray*}
Therefore it follows that
\begin{eqnarray*}
(x-y)^{[l]}*(q^my)^{[n-l]}&=&q^{m(n-l)}y^{n-l}\times \sum_{u=0}^l {l \brack u}_q (-1)^uq^{{u \choose 2}}x^{l-u}y^u\\
&=&(q^my)^{n-l}\times (x-y)^{[l]}.
\end{eqnarray*}
\qed

\begin{defi}
\label{defi:q-transform}
{\rm
(\cite{gadouleau2007, gadouleau2008})
The $q$-{\em transform} of a homogeneous polynomial\\ $a(x,y;m)=\sum_{i=0}^r a_i(m)x^{r-i}y^i$ is defined by
the homogeneous polynomial
$$
\overline{a}(x, y; m):=\sum_{i=0}^r a_i(m) y^{[i]}*x^{[r-i]}.
$$
}
\end{defi}

The following lemma is essential:

\begin{lem}
\label{lem:equations}
\begin{enumerate}
\item[{\rm (1)}] $\displaystyle
{a \brack b}_q={a-1 \brack b}_q+q^{a-b}{a-1 \brack b-1}_q=q^b{a-1 \brack b}_q+{a-1 \brack b-1}_q;
$
\item[{\rm (2)}] $\displaystyle
{a \brack b}_q {b \brack c}_q={a \brack b-c}_q {a-b+c \brack c}_q;
$
\item[{\rm (3)}] $\displaystyle
{a \brack b}_q=\sum_{i=0}^nq^{i(a-b-n+i)}{n \brack i}_q {a-n \brack b-i}_q,\:\:\:n=0,1,\ldots, a;
$
\item[{\rm (4)}] $\displaystyle
{a+b \choose 2}={a \choose 2}+ab+{b \choose 2}.
$
\end{enumerate}
\end{lem}

By combining Proposition \ref{prop:dual-enumerator} and Lemmas \ref{lem:q-products} and \ref{lem:equations}, we have the following MacWilliams type identity for a rank weight enumerator of a Delsarte rank-metric code.

\begin{thm}
\label{thm:MacWilliams}
Let ${\cal C}$ be a Delsarte rank-metric code in $\mbox{\rm Mat}(n\times m, {\mathbb F}_q)$.
Then 
$$
W_{{\cal C}^{\perp}}^{{\rm R}} (x, y)=\frac{1}{|{\cal C}|}\overline{W_{{\cal C}}^{{\rm R}}} (x+(q^m-1)y, x-y).
$$
\end{thm}
\pf
By Corollary 1 in \cite{gadouleau2008}, we have that
$$
(x-y)^{[i]}*(x+(q^m-1)y)^{[n-i]}=\sum_{j=0}^nP_j(i ; m, n)y^jx^{n-j},
$$
where
$$
P_j(i ; m, n):=\sum_{l=0}^j {l \brack j}_q {n-i \brack j-l}_q (-1)^l q^{{l \choose 2}} q^{l(n-i)} \prod_{u=0}^{j-l-1} (q^{m-l}-q^u).
$$
In addition, from Proposition \ref{prop:dual-enumerator}, it is sufficient to prove that
$$
 \sum_{l=0}^{j} {n-i \brack j-l}_q {n-j+l \brack l}_q  (-1)^l q^{{l \choose 2}} q^{m(j-l)}=P_j(i ; m, n).
$$

By using Lemma \ref{lem:equations}, we have that
\begin{eqnarray*}
\mbox{\rm L.H.S.}&=&
 \sum_{l=0}^{j} {n-i \brack j-l}_q \left\{ \sum_{u=0}^i q^{u(n-j-i+u)} {i \brack u}_q {n-j+l-i \brack l-u}_q\right\}  (-1)^l q^{{l \choose 2}} q^{m(j-l)}\\
 &=&
\sum_{u=0}^i {i \brack u}_q q^{u(n-j-i+u)}\sum_{l=0}^{j} {n-i \brack j-l}_q {n-j+l-i \brack l-u}_q  (-1)^l q^{{l \choose 2}} q^{m(j-l)}\\
&=&
\sum_{u=0}^i {i \brack u}_q q^{u(n-j-i+u)}{n-i \brack j-u}_q \sum_{l=0}^{j}  {j-u \brack l-u}_q  (-1)^l q^{{l \choose 2}} q^{m(j-l)}\\
&=&
\sum_{u=0}^i {i \brack u}_q q^{u(n-j-i+u)}{n-i \brack j-u}_q \sum_{l=0}^{j-u}  {j-u \brack l}_q  (-1)^{l+u} q^{{l+u \choose 2}} q^{m(j-l-u)}\\
&=&
\sum_{u=0}^i {i \brack u}_q q^{u(n-j-i+u)}(-1)^u q^{{u \choose 2}}{n-i \brack j-u}_q \sum_{l=0}^{j-u}  {j-u \brack l}_q  (-1)^{l} q^{{l \choose 2}} (q^{m})^{(j-u)^l} (q^u)^l\\
&=&
\sum_{u=0}^i {i \brack u}_q q^{u(n-j-i+u)}(-1)^u q^{{u \choose 2}}{n-i \brack j-u}_q \prod_{i=0}^{j-u-1} (q^m-q^{i+u})\\
&=&
\sum_{u=0}^i {i \brack u}_q q^{u(n-i)}(-1)^u q^{{u \choose 2}}{n-i \brack j-u}_q \prod_{i=0}^{j-u-1} (q^{m-u}-q^{i})\\
&=& P_j(i ; m, n).
\end{eqnarray*}
\qed

\section*{Acknowledgement}
The author would like to thank Prof. Frederique Oggier for introducing him to the concept of rank-metric codes.

\end{document}